\documentclass[leqno]{amsart}
       
\usepackage{amsmath,amsfonts,amssymb}
\usepackage{amsthm}
\usepackage{hyperref}

\def\dis
{\displaystyle}

\def\R{{\mathbb R}}
\def\N{{\mathbb N}}

\def\T{{\mathbb T}}

\def\Sch{{\mathcal S}} 
\def\virgp{\raise 2pt\hbox{,}}

\def\bu{{\bf u}}

\def\({\left(}
\def\){\right)}
\def\<{\left\langle}
\def\>{\right\rangle}
\def\le{\leqslant}
\def\ge{\geqslant}

\def\Eq#1#2{\mathop{\sim}\limits_{#1\rightarrow#2}}

\def\d{{\partial}}
\def\a{{\tt a}}

\def\eps{\varepsilon}
\def\l{\lambda}

\def\si{{\sigma}}
\def\u{{\tt u}}
\def\v{{\tt v}}

\def\F{\mathcal F}
\def\O{\mathcal O}

\DeclareMathOperator{\RE}{Re}
\DeclareMathOperator{\IM}{Im}
\DeclareMathOperator{\DIV}{div}

\theoremstyle{plain}
\newtheorem{theorem}{Theorem}[section]
\newtheorem{lemma}[theorem]{Lemma}
\newtheorem{corollary}[theorem]{Corollary}
\newtheorem{proposition}[theorem]{Proposition}
\newtheorem{hyp}[theorem]{Assumption}

\theoremstyle{definition}

\theoremstyle{remark}
\newtheorem{remark}[theorem]{Remark}
\newtheorem*{remark*}{Remark}

\numberwithin{equation}{section}

\begin{document}

\title[Semi-classical limit for NLS]{On the semi-classical limit for
  the nonlinear Schr\"odinger equation}   
\author[R. Carles]{R{\'e}mi Carles}
\address{D\'epartement de Math\'ematiques, UMR CNRS 5149\\ CC 051\\
  Universit\'e Montpellier~2\\ Place Eug\`ene Bataillon\\ 34095
  Montpellier cedex 5\\ France\footnote{Present address: Wolfgang
  Pauli Institute, Universit\"at Wien, 
        Nordbergstr.~15, A-1090 Wien}}
\email{Remi.Carles@math.cnrs.fr}
\begin{abstract}
We review some results concerning the semi-classical limit for
  the nonlinear Schr\"odinger equation, with or without an external
  potential. We consider initial data which 
  are either of the WKB type, or very concentrated as the
  semi-classical parameter goes to zero. We sketch the techniques used
  according to various frameworks, and point out some open problems. 
\end{abstract}
\maketitle

\tableofcontents
\newpage

\section{Introduction}
\label{sec:intro}

Consider the nonlinear Schr\"odinger equation (NLS):
\begin{equation}
  \label{eq:nlssemi}
  i\eps \d_t u^\eps +\frac{\eps^2}{2}\Delta u^\eps = V u^\eps
  +\eps^\kappa f\(|u^\eps|^2\)u^\eps,\quad (t,x)\in I\times \R^n,
\end{equation}
where the potential $V=V(x)$ and the
nonlinearity $f$ are 
real-valued. In some specified cases, we allow the potential to be
time-dependent. To simplify the discussion, we assume that $\kappa\ge 0$
is an integer. More precise assumptions will be made
according to the 
different cases we study. We assume $\eps\in]0,1]$, and we aim at
describing the asymptotic behavior of $u^\eps$ as $\eps \to 0$, for
the following two families of initial data:\\

\emph{Monokinetic WKB initial data:}
  \begin{equation}
    \label{eq:CIbkw}
    u^\eps(0,x) = a_0^\eps(x)e^{i\phi_0(x)/\eps}, \quad \text{with }
a_0^\eps(x)\Eq \eps 0 a_0(x) + \eps a_1(x)+\eps^2 a_2(x)+\ldots,
  \end{equation}
in the sense of asymptotic expansion. \\

\emph{Concentrated initial data:}
  \begin{equation}
    \label{eq:CIconc}
    u^\eps(0,x)= R\( \frac{x-x_0}{\eps}\)e^{i x\cdot \xi_0/\eps},
  \end{equation}
for some point $(x_0,\xi_0)$ in the phase space $\R^{2n}$,
independent of $\eps$.  
\smallbreak

There are at least two motivations for such a
study, referred to as
\emph{semi-classical analysis} or \emph{geometrical optics}. We
outline them here, and refer to the survey \cite{RauchUtah} 
for a broader discussion on this subject. The first one comes from the
applied mathematics, and may find 
its origins in physics. In the case of \eqref{eq:nlssemi}, suppose
that $\eps$ 
represents the (rescaled) Planck constant. It may be small compared to
the other parameters at stake. In this case, it is sensible to consider
that the asymptotic behavior of $u^\eps$ as $\eps\to 0$ provides a
reliable approximation of the exact solution. Hopefully, the
asymptotic model is easier to describe than the initial one
\eqref{eq:nlssemi}--\eqref{eq:CIbkw}. If $V$ is a confining potential
(e.g.\ harmonic potential), then \eqref{eq:nlssemi} may be a model to
describe Bose-Einstein condensation; see for instance
\cite{DGPS,PiSt}. The value of $\kappa$ then depends on the asymptotic
r\'egime considered. Another motivation stems from the
propagation of singularities for equations where the small parameter
$\eps$ is not necessarily present initially. Most of the studies in
this direction concern hyperbolic equations. However, 
this field is applicable to Schr\"odinger equations as well (see
e.g.\  \cite{BGT,LebeauSchrod, JeremieAIF}). The following
illustration is a straightforward consequence of the analysis
presented in \S\ref{sec:grenier}: 
\begin{theorem}[\cite{CaARMA}, Cor.~1.7]\label{theo:loss}
  Let $n\ge 3$. Consider the cubic,
  defocusing NLS:
  \begin{equation}
    \label{eq:NLScubic}
    i\d_t u +\frac{1}{2}\Delta u = |u|^2 u,\quad x\in \R^n\quad ; \quad u_{\mid
    t=0}=u_0\, .
  \end{equation}
Denote $s_c = \frac{n}{2}-1$. Let $0<s<s_c$. We can find a
  family $(u_0^\eps)_{0<\eps \le 
  1}$ in $\Sch({\R}^n)$ with 
\begin{equation*}
  \|u_0^\eps\|_{H^s({\R}^n)} \to 0 \text{ as }\eps \to 0\, ,
\end{equation*}
and $0<t^\eps \to 0$ such that the solution $u^\eps$ to
\eqref{eq:NLScubic} associated to $u_0^\eps$ satisfies: 
\begin{equation*}
  \|u^\eps(t^\eps)\|_{H^{k}({\R}^n)} \to +\infty \text{ as }\eps \to
 0\, , \ \forall k\in \left]\frac{s}{\frac{n}{2}-s}, s\right]\, .
\end{equation*}
\end{theorem}
This result was first established in \cite{CCT2} in the case
$k=s$. The fact that one can consider a broader range for $k$, in the
spirit of \cite{Lebeau05}, relies
on a fine analysis of the limit for
\eqref{eq:nlssemi}--\eqref{eq:CIbkw}, provided essentially in
\cite{Grenier98}. 
\subsection{Monokinetic WKB initial data}
\label{sec:cibkw}

In the case of initial data of the form \eqref{eq:CIbkw}, an
approximation of the form
\begin{equation}\label{eq:defBKW}
  u^\eps(t,x)\Eq \eps 0 \(\a_0(t,x)+\eps
  \a_1(t,x)+\eps^2\a_2(t,x)+\ldots\)e^{i\Phi(t,x)/\eps} 
\end{equation}
is expected. Note that only one phase and one harmonic are sought:
this is an important feature of Schr\"odinger equations with gauge
invariant nonlinearity. In the case of wave equations for instance,
the story is completely different (see e.g.\  \cite{RauchUtah} and
references therein). Note also that such an approximation must be
expected for bounded time only. Even in the linear case $f\equiv 0$, a
caustic appears in finite time in general.  Near a caustic, all the terms
$\Phi$, ${\tt a}_0$, ${\tt a}_1$, \ldots become singular. Past the
caustic, several phases are necessary in general to describe the
asymptotic behavior of the solution (see e.g.\  \cite{Du} for a general
theory in the linear case). However, we will see that
the analogous phenomenon in the nonlinear setting (say, $f(y)=y$) with
$\kappa=0$ (highly nonlinear r\'egime) might be very different. 
\smallbreak

Plug a formal expansion of the form
\eqref{eq:defBKW} into \eqref{eq:nlssemi}. Ordering the terms in
powers of $\varepsilon$, and canceling the cascade of equations thus obtained
is aimed at yielding $\Phi$, ${\tt a}_0$, ${\tt a}_1$, \ldots 
\smallbreak

Assume for a while 
that $\kappa \ge 1$. To cancel the term
of order $\O(\eps^0)$, we find 
\begin{equation*}
    \a_0\(\partial_t \Phi +\frac{1}{2}|\nabla \Phi|^2
    +V\)=0\quad ;\quad 
    \Phi_{\mid t=0}=\phi_0\, .
  \end{equation*}
Since we seek a non-trivial profile $\a_0$, we impose a stronger
condition:  $\Phi$  must
solve the \emph{eikonal equation}  
 \begin{equation*}
    \partial_t \Phi +\frac{1}{2}|\nabla \Phi|^2
    +V=0\quad ;\quad 
    \Phi_{\mid t=0}=\phi_0\, .
  \end{equation*}
Canceling the term of order
$\O(\eps^1)$, we get: 
\begin{equation*}
  \d_t {\a}_0 +\nabla \Phi\cdot \nabla {\a}_0 +
  \frac{1}{2}{\a}_0\Delta \Phi =
\left\{
  \begin{aligned}
    &0 & \text{ if }\kappa >1,\\
    &-if\left(|{\a}_0|^2\right){\a}_0& \text{ if }\kappa =1.
  \end{aligned}
\right.
\end{equation*}
We see that  the value $\kappa =1$ is critical as far as
nonlinear effects are concerned: if $\kappa>1$, no nonlinear effect is
expected at leading order, since formally, $u^\eps \sim {\tt
  a}_0 e^{i\Phi
  /\varepsilon}$, where $\Phi$ and ${\tt a}_0$ do not depend
on the nonlinearity $f$. If $\kappa =1$, then ${\tt a}_0$ solves a
nonlinear equation involving $f$. 
\smallbreak

We will see in Section~\ref{sec:weak} that when $\kappa=1$, ${\tt
  a}_0$ solves a transport 
equation that turns out to be a ordinary differential equation along
the rays of geometrical optics, as is usual in the hyperbolic case
(see e.g.\  \cite{RauchUtah}). More typical of Schr\"odinger equation is
the fact that this ordinary differential equation can be solved
explicitly.
\smallbreak 

Assume now $\kappa =0$, and proceed the same way. Plugging
\eqref{eq:defBKW} into \eqref{eq:nlssemi}, we get:
\begin{equation}
  \label{eq:bkwsurform}
\left\{
  \begin{aligned}
 \O\(\eps^0\):&\quad \d_t \Phi
 +\frac{1}{2}|\nabla\Phi|^2 + V + f\(|{\a}_0|^2\)=0,\\ 
\O\(\eps^1\):&\quad \d_t {\a}_0 +\nabla\Phi
 \cdot \nabla {\a}_0 +\frac{1}{2}{\a}_0\Delta \Phi = 
-2if'\(|\a_0|^2\) \RE\({\a}_0\overline{{\a}_1}\)\a_0. 
\end{aligned}
\right.
\end{equation}
We see that there is a strong
coupling between the phase and the main amplitude: ${\a}_0$ is present
in the equation for $\Phi$. In addition, the above system is not closed:
$\Phi$ is determined in function of ${\a}_0$, and ${\a}_0$ is
determined in function of ${\a}_1$. Even if we pursued the cascade
of equations, this phenomenon would remain: no matter how many terms
are computed, the system is never closed (see \cite{PGX93}). This is a
typical feature of supercritical cases in nonlinear geometrical optics
(see \cite{CheverryBullSMF,CG05}). We shall call the study of this
case \emph{highly nonlinear WKB analysis}.  We will see in
\S\ref{sec:strong} some ways to overcome the difficulties pointed out
above, especially in the case $f'>0$ (defocusing, cubic at the origin,
nonlinearity). 
\begin{remark}
  We consider only \emph{monokinetic} initial data. Studying the
  nonlinear effects relevant at leading order ($\kappa=0$
  or $1$) when 
  the datum is of the form 
  \begin{equation*}
    u^\eps(0,x) =
    a_0(x)e^{i\phi_0(x)/\eps}+b_0(x)e^{i\varphi_0(x)/\eps}\quad
    (\phi_0\not =\varphi_0),
  \end{equation*}
for instance, seems to be an open problem.
\end{remark}
\subsection{Concentrated initial data}
\label{sec:ciconc}
For data of the form \eqref{eq:CIconc}, a formal analysis shows that
the case $\kappa =0$ is critical: if $\kappa>0$ (not necessarily an
integer), no nonlinear effect is expected at leading order. We shall
therefore restrict our attention to the case $\kappa=0$. We also
consider the case of a pure power nonlinearity,
\begin{equation*}
  f\(|u^\eps|^2\)= \l|u^\eps|^{2\si},
\end{equation*}
for some $\si>0$ and $\l\in \R$. In this case, setting $\u^\eps =
\eps^{-n/2}u^\eps$, 
\eqref{eq:nlssemi}--\eqref{eq:CIconc} is equivalent to:
\begin{equation}
  \label{eq:nlssemiconc}
  i\eps \d_t \u^\eps +\frac{\eps^2}{2}\Delta \u^\eps = V \u^\eps
  +\l\eps^{n\si}|\u^\eps|^{2\si}\u^\eps\ ;\  \u^\eps(0,x) =
  \frac{1}{\eps^{n/2}}R\( \frac{x-x_0}{\eps}\)e^{i x\cdot
  \xi_0/\eps}. 
\end{equation}
In the case $\l>0$ (defocusing nonlinearity), dispersive effects are
expected to alter the concentrated form of the initial data. This is
proved in \cite{CaIHP,CM} when the external potential is a polynomial
of degree at most two. It seems that proving a similar result in the
more general (and fairly natural) framework of smooth, sub-quadratic,
potentials, is still an open problem. Note also that the dispersive
effect can be just the first step of the dynamics. It can be followed
by a linear dynamics induced by the potential. In this r\'egime, the
potential may cause a refocusing phenomenon. This is the case for
instance when $V$ is an isotropic harmonic potential \cite{CaIHP}. We
discuss more precisely these results in \S\ref{sec:concdefoc}.
\smallbreak

When $\l<0$ (focusing nonlinearity), several papers have considered
the case when the profile $R$ is the ground state associated to NLS
without potential, that is when $R=Q$, where $Q$ is the unique
positive, radially symmetric, solution of:
\begin{equation*}
  -\frac{1}{2}\Delta Q +Q +\l |Q|^{2\si}Q=0. 
\end{equation*}
When $\si<2/n$ (sub-critical case at the $L^2$ level), orbital stability
of the solitary wave suggests that the solution $u^\eps$ evolves under
the form
\begin{equation*}
  u^\eps(t,x) = Q\( \frac{x-x(t)}{\eps}\)e^{i x\cdot
  \xi(t)/\eps}e^{i\varphi^\eps(t)}. 
\end{equation*}
We will see that this is  the case, with $(x(t),\xi(t))$
given by the Hamiltonian flow associated to $-\frac{1}{2}\Delta+ V$:
the additional purely 
time dependent phase shift $\varphi^\eps$ is known explicitly in the
case without potential, but not in general. 
The first
mathematical result on this problem is due to J.~Bronski and
R.~Jerrard \cite{BJ00}. Refinements were then given by S.~Keraani
\cite{Keraani02,KeraaniCRAS,KeraaniAA}. We outline the approach of
\cite{KeraaniAA} in \S\ref{sec:concfoc}. Note also that the
semi-classical limit $\eps\to 0$ for \eqref{eq:nlssemiconc} is analogous to the
long time behavior for the solutions to \eqref{eq:nlssemiconc} with
$\eps=1$; see e.g.\  \cite{FGJS,JFGS}.

\section{WKB analysis for a weak nonlinearity}
\label{sec:weak}

When $\kappa\ge 1$, the first step in the WKB analysis presented in
\S\ref{sec:cibkw} consists in solving the eikonal equation. This step
relies on the Hamilton-Jacobi theory. It is well-known, at least when
the potential $V$ and the initial phase $\phi_0$ are smooth,  that the
local inversion theorem yields a local in time, smooth solution in the
neighborhood of $(t=0,x)$, for all $x\in \R^n$ (see
e.g.\  \cite{DG}). In order to have a 
local existence time which is uniform with respect to $x\in\R^n$, the
following assumption is essentially necessary (see e.g.\  \cite{CaBKW}):
\begin{hyp}\label{hyp:geom}
  The potential $V$ may depend on time: $V=V(t,x)$. We assume that the
  potential and the initial phase are smooth and sub-quadratic: 
  \begin{itemize}
  \item $V\in C^\infty(\R_t\times \R^n_x)$, and $\partial_x^\alpha V\in
  C(\R_t ;L^\infty(\R^n_x))$ as soon as $|\alpha|\ge 2$.
  \item $\phi_0\in C^\infty( \R^n)$, and $\partial_x^\alpha \phi_0\in
  L^\infty(\R^n)$ as soon as $|\alpha|\ge 2$.
  \end{itemize}
\end{hyp}
\begin{remark}
 Of course, if we worked on a compact set instead of $\R^n$, the above
assumptions would not be necessary. 
\end{remark}
A global inversion result (see \cite{SchwartzBook} or \cite{DG}) and Gronwall
lemma yield:
\begin{lemma}\label{lem:hj}
  Under Assumption~\ref{hyp:geom}, there exist $T>0$ and a unique
  solution $\phi_{\rm eik}\in C^\infty([0,T]\times\R^n)$ to:
  \begin{equation}
    \label{eq:eik}
    \d_t \phi_{\rm eik} +\frac{1}{2}|\nabla \phi_{\rm eik}|^2
    +V=0\quad ;\quad 
    \phi_{{\rm eik} \mid t=0}=\phi_0\, .
  \end{equation}
This solution is subquadratic: $\d_x^\alpha \phi_{\rm eik} \in
  L^\infty([0,T]\times\R^n)$ as soon as $|\alpha|\ge 2$. 
\end{lemma}
\begin{remark}\label{rema:sousquad}
  In \cite{CaBKW}, examples are given, that show that if either the
  potential $V$ or the initial phase $\phi_0$ has a super-quadratic
  growth at infinity, the above result fails. Sub-quadratic
  potentials play a special role in the mathematical analysis of
  Schr\"odinger equations: the results of \cite{Fujiwara79,Fujiwara}
  imply local in time Strichartz estimates for the
  semi-group associated to $-\Delta +V$. On the other hand, in space
  dimension $n=1$, $-\d_x^2-x^4$ is not essentially self-adjoint on
  $C_0^\infty(\R)$ (see \cite[Chap.~13, Sect.~6,
  Cor.~22]{Dunford}). If $V$ tends to $+\infty$ at infinity, with
  super-quadratic growth, the available results are very different
  from those of the sub-quadratic case, see e.g.\  \cite{YajZha01,YajZha04}.  
\end{remark}
To prove this lemma, we introduce the
Hamiltonian flow: 
\begin{equation}
  \label{eq:hamilton}
\left\{
  \begin{aligned}
   &\partial_t x(t,y) = \xi \left(t,y\right)\quad ;\quad x(0,y)=y,\\ 
   &\partial_t \xi(t,y) = -\nabla V\left(t,x(t,y)\right)\quad ;\quad
   \xi(0,y)=\nabla \phi_0(y).
  \end{aligned}
\right.
\end{equation}
The time $T$ is such that the map $y\mapsto x(t,y)$ is a
diffeomorphism of $\R^n$ for $t\in [0,T]$. 
Therefore, the Jacobi determinant
\begin{equation*}
  J_t(y) ={\rm det}\nabla_y x(t,y),
\end{equation*}
is bounded from above, and from below away from zero, for
$t\in[0,T]$. 
The justification of the leading order asymptotics sketched in
\S\ref{sec:cibkw} is:
\begin{proposition}\label{prop:sub}
  Let
  $\kappa \ge 1$ and $f\in C^\infty(\R_+;\R)$. Assume that there
  exists a smooth function $a_0$ 
  independent of $\eps$ such that
\begin{equation*}
  a_0^\eps \to a_0 \text{ in }H^s(\R^n),\quad \forall s\ge 0. 
\end{equation*}
Then under Assumption~\ref{hyp:geom}, for all $\eps \in ]0,1]$,
  \eqref{eq:nlssemi}--\eqref{eq:CIbkw} has a unique solution $u^\eps \in
  C^\infty([0,T]\times \R^n)\cap C([0,T];H^s)$ for all $s>n/2$, where
  $T$ is given by Lemma~\ref{lem:hj}. Moreover,
  there exist $ a, G\in C^\infty([0,T]\times \R^n)$,
  independent of $\eps \in ]0,1]$, where
  $ a \in C([0,T]; L^2\cap L^\infty)$, and $G$ is
  real-valued with $G \in C([0,T];   L^\infty)$, such that:
  \begin{equation*}
    \left\| u^\eps -  a e^{i\eps^{\kappa -1}G}e^{i\phi_{\rm eik}
        /\eps}\right\|_{L^\infty([0,T]; L^2\cap L^\infty) } \to 0\quad
    \text{as }\eps \to 0.
  \end{equation*}
The profile $ a$ solves the initial value problem:
\begin{equation}\label{eq:alibre}
  \partial_t  a +\nabla \phi_{\rm eik}\cdot \nabla a +
  \frac{1}{2} a\Delta \phi_{\rm eik} =0\quad ;\quad a_{\mid
    t=0}=a_0, 
\end{equation}
and $G$  depends nonlinearly on $ a$:
\begin{equation*}
  \begin{aligned}
    a(t,x) &= \frac{1}{\sqrt{J_t(y(t,x))}}a_0\left(y(t,x)\right),\\
  G(t,x) &= -\int_0^t
  f\( J_s(y(t,x))^{-1}\left|a_0(y(t,x)) \right|^2\) ds.
  \end{aligned}
\end{equation*}
In particular, if $\kappa>1$, then 
\begin{equation*}
    \left\| u^\eps - a e^{i\phi_{\rm eik}
        /\eps}\right\|_{L^\infty([0,T]; L^2\cap L^\infty) } \to 0\quad
    \text{as }\eps \to 0,
  \end{equation*}
and no nonlinear effect is present in the leading order behavior of
$u^\eps$. If $\kappa =1$, nonlinear effects are present at leading
order, measured by $G$.
\end{proposition}
We see that the
critical nonlinear effect (case $\kappa=1$) is a self-modulation of
the amplitude. In 
the context of laser physics, this phenomenon is known as
\emph{phase self-modulation} (see e.g.\  \cite{ZS,Boyd,Donnat}). 
\begin{proof}[Sketch of the proof]
  The proof given in \cite{CaBKW} consists in changing the unknown
  function, by setting
  \begin{equation*}
    a^\eps = u^\eps e^{-i\phi_{\rm eik}/\eps},
  \end{equation*}
where $\phi_{\rm eik}$ is given by Lemma~\ref{lem:hj}. Then
\eqref{eq:nlssemi}--\eqref{eq:CIbkw} is equivalent to:
\begin{equation*}
    \d_t a^\eps +\nabla \phi_{\rm eik}\cdot \nabla
  a^\eps +\frac{1}{2}a^\eps\Delta \phi_{\rm eik}
  =i\frac{\eps}{2}\Delta a^\eps 
  -i\eps^{\kappa -1}f\(|a^\eps|^2\) a^\eps\quad ;\quad 
a^\eps_{\mid  t=0}=a_0^\eps. 
\end{equation*}
Energy estimates show that the above equation has a unique, smooth
solution $a^\eps\in C([0,T];H^s)$ for all $s>n/2$, uniformly bounded
for $\eps\in ]0,1]$. This step uses the facts that $\phi_{\rm eik}$ is
sub-quadratic and $i\Delta$ is skew-symmetric. We can then neglect the
terms $\eps \Delta a^\eps$ and $a_0^\eps -a_0$, so that $\|a^\eps-
\widetilde a^\eps\|_{L^\infty([0,T];H^s)}=o(1)$, where: 
\begin{equation}\label{eq:atildeeps}
  \d_t \widetilde a^\eps +\nabla \phi_{\rm eik}\cdot
  \nabla \widetilde a^\eps + 
  \frac{1}{2}\widetilde a^\eps\Delta \phi_{\rm eik} =
  -i\eps^{\kappa -1}f\(|\widetilde
  a^\eps|^2\)\widetilde a^\eps \quad ; \quad 
\widetilde a^\eps_{\mid  t=0}=a_0.
\end{equation}
Recall that $J_t(y)$ is the Jacobi determinant. Denote
\begin{equation*}
  A^\eps(t,y) := \widetilde a^\eps \left(t, x(t,y)
  \right)\sqrt{J_t(y)}. 
\end{equation*}
We see that so long as $y\mapsto x(t,y)$ defines a global
diffeomorphism (which is guaranteed for $t\in [0,T]$ by construction),
\eqref{eq:atildeeps} is equivalent to:  
\begin{equation*}
  \partial_t A^\eps = -i\eps^{\kappa
  -1}f\(J_t(y)^{-1}\left|A^\eps \right|^2\)
  A^\eps\quad ; \quad A^\eps(0,y)=a_0(y). 
\end{equation*}
This ordinary differential equation along the rays of geometrical
optics can be solved explicitly, after we have remarked the identity
$\d_t |A^\eps|^2=0$: 
\begin{equation*}
  A^\eps(t,y) = a_0(y) \exp\left(-i\eps^{\kappa -1} \int_0^t
  f\(J_s(y)^{-1}\left|a_0(y) \right|^2\) ds\right). 
\end{equation*}
Back to the initial solution $u^\eps$, this yields the proposition.
\end{proof}
\begin{remark}
  A similar result is proved in \cite{CaMaSp} for the equation
  \begin{equation*}
   i\eps \d_t u^\eps +\frac{\eps^2}{2}\Delta u^\eps = V(x) u^\eps
  +V_\Gamma\(\frac{x}{\eps}\)u^\eps+\lambda \eps |u^\eps|^{2\si}u^\eps,
  \end{equation*}
where $V_\Gamma$ is lattice-periodic. The presence of this rapidly
oscillatory potential changes dramatically the geometry of the
propagation. Using the corresponding Bloch theory, a similar phase
self-modulation phenomenon is proved, under the assumption that the
initial data are well-prepared. Removing this assumption, or
considering highly nonlinear r\'egimes (as in \S\ref{sec:strong}) are
interesting open questions, and have physical motivations in the
context of Bose--Einstein condensation. 
\end{remark}

\section{Highly nonlinear WKB analysis: $\kappa=0$}
\label{sec:strong}

We saw in \S\ref{sec:cibkw} that constructing a formal asymptotic
expansion for \eqref{eq:nlssemi}--\eqref{eq:CIbkw} is a delicate issue
when $\kappa=0$. We also point out that another problem arises,
even if one has managed to construct an approximate solution $v^\eps$ that
solves
\begin{equation}\label{eq:approx}
  i\eps \d_t v^\eps +\frac{\eps^2}{2}\Delta v^\eps = V v^\eps
  +f\(|v^\eps|^2\)v^\eps+\eps^N r_N^\eps \quad ;\quad
  v^\eps_{\mid t=0}=u^\eps_{\mid t=0}, 
\end{equation}
where $N$ is large, and $r_N^\eps$ is bounded in $L^2$ for
instance. Setting $w^\eps = u^\eps -v^\eps$, and supposing that
$u^\eps$ and $v^\eps$ remain bounded in $L^\infty(\R^n)$ on a time
interval $[0,t]$, the usual $L^2$ estimate for Schr\"odinger equations
yields:
\begin{equation*}
  \eps\|w^\eps(t)\|_{L^2} \le C \int_0^t \|w^\eps(\tau)\|_{L^2}d\tau +
  2 \eps^N \int_0^t \|r_N^\eps(\tau)\|_{L^2}d\tau.
\end{equation*}
We infer, using Gronwall lemma:
\begin{equation*}
  \|w^\eps(t)\|_{L^2}\le C \eps^{N-1}e^{Ct/\eps}.
\end{equation*}
The exponential factor shows that this method may yield interesting
results only up to time of the order $c\eps |\log \eps|^\theta$ for
some $c,\theta>0$. Note that in some functional analysis contexts, this
may be satisfactory  (see e.g.\  \cite{CCT2}, or the appendices in
\cite{BGTENS,CaARMA}). However, it seems reasonable to wish to have a
description of the solution of \eqref{eq:nlssemi}--\eqref{eq:CIbkw} at
least on a time interval independent of $\eps$. We list below several
approaches that yield such information. 
\begin{remark}\label{rema:instab}
  In a slightly different context, a fairly explicit example 
  in \cite{CaJHDE} shows that 
  one may find a function satisfying \eqref{eq:approx} for $N$
  arbitrarily large, such that
  \begin{equation*}
    \liminf_{\eps \to
    0}\|u^\eps(t^\eps)-v^\eps(t^\eps)\|_{L^2}>0,\quad \text{for
    }t^\eps=\eps^\beta\text{ and }0<\beta<1. 
  \end{equation*}
Therefore, the stability
issue in this highly nonlinear r\'egime is really delicate.
\end{remark}
\subsection{Modulated energy functional}
\label{sec:brenier}

A general technique was introduced by Y.~Brenier in
\cite{BrenierCPDE}. It yields the convergence of some physically
important quantities (such as the Wigner measure, see
e.g.\  \cite{GMMP,LionsPaul}), but not of the wave function $u^\eps$ itself. In
the case of the nonlinear Schr\"odinger equation, it has been used by
P.~Zhang \cite{ZhangJPDE} (see also \cite{ZhangSIMA} for the case of
the Schr\"odinger--Poisson equation). More recently, F.~Lin and
P.~Zhang have adapted this approach in the case of the
Gross-Pitaevskii equation, in the exterior of an obstacle
\cite{LinZhang}. We shall present the technique of Brenier in the case
of \eqref{eq:nlssemi}--\eqref{eq:CIbkw}, using the simplified approach
of  \cite{LinZhang}. In all this paragraph, we will assume $V\equiv
0$: no external potential is present.
\smallbreak

The first step consists in guessing a suitable approximate
solution. Even though the system \eqref{eq:bkwsurform} is not closed,
the analysis of \S\ref{sec:weak} shows that so long as $\Phi$ is
smooth and $\nabla \Phi$ is a global diffeomorphism, the second
equation of \eqref{eq:bkwsurform} is of the form:
\begin{equation*}
  \dot \a_0 = i \Xi\a_0,
\end{equation*}
where $\dot \a_0$ stands for the differentiation along the rays
associated to $\nabla \Phi$, and $\Xi$ is real-valued. In particular,
the modulus of $\a_0$ is constant along these rays. Setting
$(\rho,v)=(|\a_0|^2,\nabla \Phi)$ as a new unknown function,
\eqref{eq:bkwsurform} yields:
\begin{equation}\label{eq:euler}
\left\{
  \begin{aligned}
    &\d_t v+v\cdot \nabla v + f'(\rho)\nabla \rho =0\quad ;\quad
    v_{\mid t=0}=\nabla\phi_0,\\
    &\d_t \rho +v\cdot \nabla \rho +\rho \DIV v =0\quad ;\quad
    \rho_{\mid t =0}=|a_0|^2. 
  \end{aligned}
\right.
\end{equation}
If $f'>0$, we get a compressible Euler equation, which  is hyperbolic
symmetric in the sense of Friedrichs. We shall assume now that
$f'\equiv 1$, that is, we consider a cubic, defocusing nonlinearity in
\eqref{eq:nlssemi}. Note that older formal approaches suggest the
introduction of \eqref{eq:euler} as a limiting equation. In
\cite[Chap.~III]{LandauQ}, we find:
\begin{equation}\label{eq:landau}
  \left\{
\begin{aligned}
    &\d_t \Phi^\eps +\frac{1}{2}\left|\nabla
    \Phi^\eps\right|^2 + 
    |\a^\eps|^2= \eps^2 \frac{\Delta
    \a^\eps}{2 \a^\eps}\quad ; \quad 
    \Phi^\eps_{\mid t=0}=\phi_0\, ,\\
&\d_t \a^\eps +\nabla \Phi^\eps \cdot \nabla
    \a^\eps +\frac{1}{2}\a^\eps 
\Delta \Phi^\eps  = 0\quad  ;\quad
\a^\eps_{\mid t=0}= a^\eps_0\, . 
\end{aligned}
\right.
\end{equation}
Of course, this choice is not adapted when the amplitude $\a^\eps$
vanishes, so it must be left out for a rigorous mathematical analysis,
when $a^\eps_0 \in L^2(\R^n)$. Passing formally to the limit $\eps\to
0$, the right hand side of the equation for $\Phi^\eps$ vanishes, and
using the hydrodynamical variables as above, we retrieve
\eqref{eq:euler}. 
\smallbreak

The modulated energy functional associated to
\eqref{eq:nlssemi}--\eqref{eq:CIbkw} when $V\equiv 0$ and $f(y)=y$ is:
\begin{equation*}
  H^\eps(t)= \frac{1}{2}\int_{\R^n} \left| (\eps \nabla
  -iv)u^\eps(t,x)\right|^2dx +\frac{1}{2}\int_{\R^n} \(\rho^\eps(t,x)
  -\rho(t,x)\)^2dx, 
\end{equation*}
where we have set $\rho^\eps = |u^\eps|^2$. We find that the time
derivative of this  modulated energy functional is:
\begin{equation*}
  \begin{aligned}
    \frac{d}{dt}H^\eps(t) =& \frac{\eps^2}{4}\int \nabla (\DIV v)\cdot \nabla
  \rho^\eps - \sum_{j,k} \int \d_j v_k \RE\( (\eps \d_j -iv_j)u^\eps
  \overline{(\eps \d_k -iv_k)u^\eps}\)\\
& + \frac{3}{2}\int \(\rho^\eps
  -\rho\)^2 \DIV  v  .  
  \end{aligned}
\end{equation*}
The last two terms are estimated by $\|\nabla
v(t)\|_{L^\infty}H^\eps(t)$. For the first term, write
\begin{align*}
  \eps^2\int \DIV \(\nabla v\)\cdot \nabla |u^\eps|^2 &=  \eps\int
  \DIV \(\nabla v\)\cdot \(\overline u^\eps \eps\nabla u^\eps +
  u^\eps \eps\nabla\overline u^\eps \)\\
&= \eps\int
  \DIV \(\nabla v\)\cdot\( \overline u^\eps (\eps\nabla-iv)u^\eps
  +u^\eps \overline {(\eps\nabla-iv)u^\eps} \).  
\end{align*}
Since $\|u^\eps(t)\|_{L^2}=\|a_0\|_{L^2}$  and
$v\in L^\infty([0,T];W^{2,\infty})$, Young's inequality yields:
\begin{equation*}
  \frac{d}{dt}H^\eps(t) \le C\(H^\eps(t) +\eps^2\),
\end{equation*}
so long as $v$ remains smooth, that is, before shocks appear in
\eqref{eq:euler}. We conclude thanks to Gronwall lemma:
\begin{theorem}\label{theo:mef}
  Let $n\ge 1$, and assume that $\kappa=0$, $V\equiv 0$ and
  $f(y)=y$. Assume that there 
  exists a smooth function $a_0$ 
  independent of $\eps$ such that
\begin{equation*}
  a_0^\eps \to a_0 \text{ in }H^s(\R^n),\quad \forall s\ge 0. 
\end{equation*}
Assume also that $\phi_0\in C^\infty(\R^n;\R)$ is such that $\nabla
\phi_0\in H^s(\R^n)$ for all $s\ge 0$. Then there exists $T>0$
independent of $\eps>0$ such that \eqref{eq:nlssemi}--\eqref{eq:CIbkw}
has a unique solution $u^\eps \in 
  C^\infty([0,T]\times \R^n)\cap C([0,T];H^s)$ for all $s>n/2$. In
  addition, as $\eps \to 0$, 
  \begin{equation*}
    \left\| (\eps \nabla
  -iv)u^\eps\right\|_{L^\infty([0,T];L^2)}^2 + \left\| |u^\eps|^2
  -\rho\right\|_{L^\infty([0,T];L^2)}^2  =\O\(\eps^2+
  \left\||a_0^\eps|^2 - |a_0|^2\right\|_{L^2}^2\).  
  \end{equation*}
\end{theorem}
In the above theorem, we have not tried to compute the lowest possible
value for the Sobolev regularity $s$ given by the proof, nor  shall we
try in the other sections. 
\begin{remark}
  In the more general case where the nonlinearity is 
  $f(y)=y^\si$, with $\si\in\N$, a generalization of the above
  modulated energy functional was 
  introduced in \cite{AC-P}. In particular, the analogue of
  Theorem~\ref{theo:mef} is proved. This includes for instance the
  case of the quintic, defocusing nonlinearity.
\end{remark}

One might be afraid that the above result is somehow contradictory with
Remark~\ref{rema:instab}, or with the results of \cite{CaARMA}. A
typical example in \cite{CaARMA}, under the assumptions of
Theorem~\ref{theo:mef}, consists in choosing $a_0^\eps=a_0$
independent of $\eps$, and considering $v^\eps$ solving
\eqref{eq:nlssemi}--\eqref{eq:CIbkw} with $\widetilde a_0^\eps =
(1+\eps^{1-\alpha})a_0$ ($0<\alpha<1$). Then for $t^\eps$ of order
$\eps^\alpha$,  
\begin{equation*}
  \liminf_{\eps \to
  0}\left\|u^\eps(t^\eps)-v^\eps(t^\eps)\right\|_{L^2}>0. 
\end{equation*}
Yet, there is no contradiction with Theorem~\ref{theo:mef}: the
instability mechanism in \cite{CaARMA} is the appearance of an extra
oscillatory factor in $v^\eps$. This oscillation shows up  essentially
through a multiplicative factor of the form $e^{ig(t,x)/\eps^\alpha}$.
It does not affects the 
modulus of the wave function, and vanishes in the limit $\eps \to 0$
of $\eps\nabla v^\eps$. \smallbreak

We can therefore conclude that the modulated
energy functional shares several features with the Wigner measure. It
is a rather general tool: in \cite{LinZhang}, the authors consider a
nonlinear Schr\"odinger equation with a boundary condition, aspect
which apparently cannot be recovered with the approach of E.~Grenier
recalled in the next paragraph. On the other hand, by definition, it
ignores the oscillatory phenomena that occur at a scale of order
$\eps^\alpha$ for $0<\alpha<1$ (for instance). The next section shows
how to get a more precise description, under similar assumptions.

\subsection{Point-wise asymptotics without potential}
\label{sec:grenier}
In this paragraph, we keep assuming $V\equiv 0$. Note that in
\eqref{eq:landau}, the supposedly small term on the right hand side is
of order $\eps^2$, while $\eps$ should be enough to neglect a term in
the limit $\eps\to 0$. We have seen in \S\ref{sec:cibkw} that the
equation for the phase is obtained after simplification by the leading
order amplitude. This explains the singular factor on the right hand
side of \eqref{eq:landau}. The main technical ingredient in
\cite{Grenier98} consists in shifting the source term in
\eqref{eq:landau} to the next order, that is, the equation for the
amplitude: we now seek $u^\eps =a^\eps e^{i\Phi^\eps/\eps}$, where the
amplitude $a^\eps$ is \emph{complex-valued} (even if 
$a_0^\eps$ is real-valued), $\Phi^\eps$ is real-valued, and:
\begin{equation}\label{eq:systemmanuel}
\left\{
  \begin{aligned}
    &\d_t \Phi^\eps +\frac{1}{2}\left|\nabla
    \Phi^\eps\right|^2 + f\left( 
    |a^\eps|^2\right)= 0\quad ; \quad
    \Phi^\eps_{\mid t=0}=\phi_0\, ,\\
&\d_t a^\eps +\nabla \Phi^\eps \cdot \nabla
    a^\eps +\frac{1}{2}a^\eps 
\Delta \Phi^\eps  = i\frac{\eps}{2}\Delta
    a^\eps\quad  ;\quad a^\eps_{\mid t=0}= a^\eps_0\, . 
  \end{aligned}
\right.
\end{equation}
Another originality of this approach lies in the fact that the phase
$\Phi^\eps$ depends on $\eps$, through the coupling of
the two equations. The idea of E.~Grenier consists in somehow
performing the usual WKB analysis ``the other 
way round'': first, solve \eqref{eq:systemmanuel}, then show that
$\Phi^\eps$ and $a^\eps$ have asymptotic expansions as $\eps \to
0$. In particular, this resolves the stability issue pointed out at the
beginning of \S\ref{sec:strong}. 
\smallbreak

To solve \eqref{eq:systemmanuel}, consider the new unknown 
\begin{equation*}
  \bu^\varepsilon = \left(
    \begin{array}[l]{c}
       \RE a^\eps \\
       \IM a^\eps \\
       \nabla \Phi^\eps 
    \end{array}
\right)\in \R^{n+2}.
\end{equation*} 
The system \eqref{eq:systemmanuel} is equivalent to a quasi-linear
equation of the form:
\begin{equation}
  \label{eq:systhyp}
  \d_t \bu^\eps +\sum_{j=1}^n
  A_j(\bu^\eps)\d_j \bu^\eps
  = \frac{\eps}{2} L 
  \bu^\eps\, , \quad\text{with}\quad L = \left(
  \begin{array}[l]{ccccc}
   0  &-\Delta &0& \dots & 0   \\
   \Delta  & 0 &0& \dots & 0  \\
   0& 0 &&0_{n\times n}& \\
   \end{array}
\right),
\end{equation}
\begin{equation*}
 \text{and}\quad A(\bu,\xi)=\sum_{j=1}^n A_j(\bu)\xi_j
= \left(
    \begin{array}[l]{ccc}
      v\cdot \xi & 0& \frac{1}{2}\RE a \,^{t}\xi \\ 
     0 &  v\cdot \xi & \frac{1}{2}\IM a\,^{t}\xi \\ 
     2f'  \RE a \, \xi
     &2f'  \IM a\, \xi &  v\cdot \xi I_n 
    \end{array}
\right),
\end{equation*}
where $f'$ stands for $f'(|a|^2)$. The system \eqref{eq:systhyp} is
hyperbolic symmetric when $f'>0$, and we can consider the following
symmetrizer:
\begin{equation*}
  S=\left(
    \begin{array}[l]{cc}
     I_2 & 0\\
     0& \frac{1}{4f'(|a|^2)}I_n
    \end{array}
\right),
\end{equation*}
which is symmetric and positive for $f'>0$. 
\begin{remark}
The argument of
$f'$ is morally bounded (this will result from the analysis), but may
have zeroes: the assumption $f'\ge 0$ cannot be considered by this
approach. For instance, justifying a WKB analysis for the quintic,
defocusing NLS remains an open problem. 
\end{remark}
An advantage for this choice of $S$ is that $SL$ remains a
skew-symmetric operator: the possible loss of derivative caused by the
second order operator $L$ does not affect the usual energy estimates
in $H^s(\R^n)$. One can then prove existence and uniqueness for
\eqref{eq:systhyp} in Sobolev spaces of sufficiently large
order. Since the last $n$ components define initially, and remain,  an
irrotational function, this 
implies that we can solve \eqref{eq:systemmanuel}. The natural limit
is given by:
\begin{equation}\label{eq:systemmanuellimite}
\left\{
  \begin{aligned}
    &\d_t \Phi +\frac{1}{2}\left|\nabla
    \Phi\right|^2 + f\left( 
    |a|^2\right)= 0\quad ; \quad
    \Phi_{\mid t=0}=\phi_0\, ,\\
&\d_t a +\nabla \Phi \cdot \nabla
    a +\frac{1}{2}a
\Delta \Phi  = 0\quad  ;\quad a_{\mid t=0}= a_0\, . 
  \end{aligned}
\right.
\end{equation}
Local existence in Sobolev spaces for \eqref{eq:systemmanuellimite}
follows from the same arguments, and one has:
\begin{theorem}\label{theo:emmanuel}
  Let $n\ge 1$, and assume that $\kappa=0$, $V\equiv 0$ and
  $f\in C^\infty(\R_+;\R)$ with $f'>0$. Assume that there 
  exists a smooth function $a_0$ 
  independent of $\eps$ such that
\begin{equation*}
  a_0^\eps \to a_0 \text{ in }H^s(\R^n),\quad \forall s\ge 0. 
\end{equation*}
Assume also that $\phi_0\in C^\infty(\R^n;\R)$ is such that $\nabla
\phi_0\in H^s(\R^n)$ for all $s\ge 0$. Then there exists $T>0$ 
independent of $\eps>0$ such that \eqref{eq:nlssemi}--\eqref{eq:CIbkw}
has a unique solution $u^\eps= a^\eps e^{i\Phi^\eps/\eps}$ in  
 $C^\infty([0,T]\times \R^n)\cap C([0,T];H^s)$ for all
 $s>n/2+2$. Moreover, $a^\eps$ and  
$\Phi^\eps$ are bounded in $L^\infty([0,T];H^s)$,
uniformly in $\eps\in ]0,1]$ and, for all $s>n/2+1$, there exists
$C_s$ such that  
  \begin{equation*}
    \left\| \nabla \Phi^\eps -\nabla \Phi\right\|_{L^\infty
    ([0,T];H^s)} +\left\| a^\eps -a\right\|_{L^\infty
    ([0,T];H^s)} \le C_s\( \eps +\left\| a^\eps_0 -a_0\right\|_{H^s}\). 
  \end{equation*}
Therefore, $\dis \left\| \Phi^\eps(t) -\Phi(t)\right\|_{H^s}  \le \widetilde
C_s t\( \eps +\left\| a^\eps_0 
    -a_0\right\|_{H^s}\)$, $\forall t\in [0,T]$.
\end{theorem}
Theorem~\ref{theo:emmanuel} does not suffice to describe the
asymptotic behavior of $u^\eps$ on the time interval $[0,T]$ though:
\begin{align*}
  u^\eps - a e^{i\Phi/\eps} =a^\eps e^{i\Phi^\eps /\eps}- a
  e^{i\Phi/\eps}
 = \(a^\eps -a\) e^{i\Phi^\eps /\eps} + a\( e^{i\Phi^\eps /\eps}-
  e^{i\Phi/\eps} \). 
\end{align*}
Therefore, we have
\begin{align*}
  \left| u^\eps - a e^{i\Phi/\eps}\right| \le |a^\eps -a| + 2|a|\left|
  \sin\( \frac{\Phi^\eps-\Phi}{2\eps}\)\right|
\end{align*}
Taking the $L^2$ norm, we infer:
\begin{align*}
 \left\| u^\eps(t) - a(t) e^{i\Phi(t)/\eps}\right\|_{L^2}& \le \|a^\eps(t)
  -a(t)\|_{L^2} + 2\|a(t)\|_{L^2}\left\| 
  \sin\( \frac{\Phi^\eps(t)-\Phi(t)}{2\eps}\)\right\|_{L^\infty}\\
& \lesssim \(\eps + \left\| a^\eps_0
 -a_0\right\|_{H^s}\)\(1+\frac{t}{\eps}\),
\end{align*}
for $s>n/2+1$. Even if $a^\eps_0
 -a_0 =\O(\eps^N)$ for $N$ large, the above estimate shows that $a
 e^{i\Phi/\eps}$ is a good approximation of $u^\eps$ as $t\to
 0$, but not necessarily at time $t=T$ for instance. To have a better
 error estimate, it is necessary to compute the next 
term in the asymptotic expansion of $(\phi^\eps,a^\eps)$ in powers of
$\eps$. Assume furthermore that
there exists $a_1\in \cap_{s\ge 0}H^s$ such that 
 \begin{equation}\label{eq:ciprec}
    a_0^\eps = a_0 +\eps a_1 +o(\eps)\quad
    \text{in }H^s, \ \forall s\ge 0.
  \end{equation}
For times of order $\O(1)$, the initial corrector $a_1$ must
be taken into account: 
\begin{proposition}\label{prop:correc}
  Define
  $(a^{(1)},\Phi^{(1)})$ by 
\begin{equation*}
\left\{
  \begin{aligned}
    &\d_t \Phi^{(1)} +\nabla \Phi \cdot \nabla \Phi^{(1)} +
    2f'\(|a|^2\)\RE\left(\overline a a^{(1)}\right)=0 , \\ 
   &\d_t a^{(1)} +\nabla\Phi\cdot \nabla a^{(1)} + \nabla
   \Phi^{(1)}\cdot \nabla a + \frac{1}{2} a^{(1)}\Delta \Phi
   +\frac{1}{2}a\Delta \Phi^{(1)}      = \frac{i}{2}\Delta a,\\
&    \Phi^{(1)}_{\mid t=0}=0\quad ;
    \quad a^{(1)}_{\mid t=0}=a_1. 
  \end{aligned}
\right.
\end{equation*}
Then $a^{(1)},\Phi^{(1)}\in
L^\infty([0,T];H^s)$ for every $s\ge 0$, and
\begin{equation*}
  \|a^\eps - a - \eps a^{(1)}\|_{L^\infty([0,T];H^s)}+
  \|\Phi^\eps - \Phi - \eps 
  \Phi^{(1)}\|_{L^\infty([0,T];H^s)} \le C_s\eps^2,\quad
  \forall s\ge 0\, . 
\end{equation*}
\end{proposition}
Despite 
the notations, it seems unadapted to consider $\Phi^{(1)}$ as being
part of the phase. Indeed, we infer from Proposition~\ref{prop:correc}
that 
\begin{equation*}
  \left\|u^\eps - a e^{i\Phi^{(1)}}
    e^{i\Phi/\eps}\right\|_{L^\infty([0,T];L^2\cap
    L^\infty)}= \O(\eps).  
\end{equation*}
Relating this information to the WKB methods presented in
\S\ref{sec:cibkw}, we would have:  
\begin{equation*}
  {\a}_0 = a e^{i\Phi^{(1)}}.
\end{equation*}
Since $\Phi^{(1)}$ depends on $a_1$ while $a$ does not, we retrieve
the fact that in super-critical r\'egimes, the leading order amplitude
in WKB methods depends on the initial first corrector $a_1$. 
\begin{remark}
  The term $e^{i\Phi^{(1)}}$ does not appear in the Wigner measure
  of $a e^{i\Phi^{(1)}} e^{i\Phi/\eps}$. Thus, from the point of view of
    Wigner measures, the asymptotic behavior of the exact solution is
    described by the Euler-type system \eqref{eq:euler}. 
\end{remark}
\begin{remark}
  If we assume that $a_0$ is real-valued, then so is $a$. If moreover
  $a_1$ is purely imaginary (for instance, if $a_1=0$), then we see
  that $a^{(1)}$ is purely imaginary, hence, $\Phi^{(1)}\equiv
  0$. 
\end{remark}
\begin{remark}
  The proof of Theorem~\ref{theo:loss} follows. Consider
  initial data of the form
  \begin{equation*}
    u_0(x) =
    \lambda^{-\frac{n}{2}+s}a_0\left(\frac{x}{\lambda}\right), \quad
    \lambda \to 0 .
  \end{equation*}
Set $\eps = \lambda^{s_c-s}$: $\eps$ and $\lambda$ go
simultaneously to zero, by assumption. Define 
\begin{equation*}
  \psi^\eps(t,x) = \lambda^{\frac{n}{2}-s}u\left(
  \lambda^{\frac{n}{2}+1-s} t,\lambda x\right)\, . 
\end{equation*}
It solves:
\begin{equation}\label{eq:psi52}
  i\eps\d_t \psi^\eps +\frac{\eps^2}{2}\Delta \psi^\eps =
  |\psi^\eps|^{2}\psi^\eps \quad 
  ;\quad \psi^\eps_{\mid t=0} = a_0(x) \, .
\end{equation}
The idea of the proof is that for times of order ${\mathcal O}(1)$,
$\psi^\eps$ has become $\eps$-oscillatory. This is rather clear from
\eqref{eq:systemmanuellimite}: even though $\Phi_{\mid t=0}=0$, we
have $\d_t \Phi_{\mid t=0}\not=0$, and rapid oscillations at scale
$\eps$ appear instantly. Back to $u$, this yields the theorem (up to
replacing $a_0$ by $|\log \lambda|^{-1}a_0$). 
\end{remark}
To conclude this paragraph, we point out an open problem concerning
the time $T_c$ when shocks appear for \eqref{eq:euler}. First, the break-up
for \eqref{eq:euler} does not allow us to deduce anything concerning
the behavior of the solution of \eqref{eq:systemmanuel}. More
generally, the notion of caustic in this case is not so
clear. Geometrically, as $t\to T_c$, the rays
for  \eqref{eq:systemmanuellimite} tend to form an envelope. In the
linear case $f\equiv 0$, this geometrical phenomenon goes along with
an analytical one: 
\begin{equation*}
  \liminf_{\eps \to 0}\|u^\eps(t)\|_{L^\infty} \to +\infty\quad
  \text{as }t\to T_c.
\end{equation*}
For instance, $\|u^\eps(t)\|_{L^\infty}\thickapprox (\eps+
|T_c-t|)^{-n/2}$ for all $t$ in
the case of a focal point (all the rays meet at one point as $t\to
T_c$). 
\smallbreak

It is not clear at all that a similar phenomenon occurs for
\eqref{eq:nlssemi} when $\kappa =0$. Suppose for instance that the
nonlinearity is cubic, defocusing, $f(y)=y$, and that the initial
profile $a_0^\eps$ does not depend on $\eps$, $a_0^\eps =a_0$. The
standard conservations of mass and energy for nonlinear Schr\"odinger
equations yield:
\begin{align*}
  &\|u^\eps(t)\|_{L^2} = \|a_0\|_{L^2}=\O(1),\\
& \|\eps\nabla u^\eps(t)\|_{L^2}^2 + \|u^\eps(t)\|_{L^4}^4
=\|\eps\nabla a_0+ia_0\nabla \phi_0\|_{L^2}^2 + \|a_0\|_{L^4}^4=\O(1). 
\end{align*}
In space dimension $n\le 3$, the solution $u^\eps$ remains in
$H^1(\R^n)$ for all time, therefore we know that the $L^2$ and $L^4$
norms of $u^\eps(t,\cdot)$ remain bounded by a constant independent of
$\eps$. This \emph{suggests} that the $L^\infty$ norm of
$u^\eps(t,\cdot)$ may remain bounded, if we can somehow inverse the
H\"older inequality
\begin{equation*}
  \|u^\eps(t)\|_{L^4}^4\le \|u^\eps(t)\|_{L^2}^2
  \|u^\eps(t)\|_{L^\infty}^2.  
\end{equation*}
One could then distinguish two
notions of caustic: a geometrical one (present in all the cases), and
an analytical one (possibly absent in the highly nonlinear case).  
\subsection{Point-wise asymptotics with an external potential}
\label{sec:cabkw}

Physical motivations may lead to the study of
\eqref{eq:nlssemi}--\eqref{eq:CIbkw} when the external potential $V$
is not zero. Mathematically, a special role is played by sub-quadratic
potentials, as we have noticed in \S\ref{sec:weak}; see
Remark~\ref{rema:sousquad}. We therefore suppose that
Assumption~\ref{hyp:geom} is satisfied. 
\smallbreak

The analysis presented in \S\ref{sec:brenier} suggests that in this
case, we have to consider solutions to a compressible Euler equation
with (possibly) unbounded external force and initial velocity:
\begin{equation}\label{eq:eulerV}
\left\{
  \begin{aligned}
    &\d_t v+v\cdot \nabla v + \nabla V+ f'(\rho)\nabla \rho =0\quad ;\quad
    v_{\mid t=0}=\nabla\phi_0,\\
    &\d_t \rho +v\cdot \nabla \rho +\rho \DIV v =0\quad ;\quad
    \rho_{\mid t =0}=|a_0|^2. 
  \end{aligned}
\right.
\end{equation}
The existence of such solutions is not standard. The
naive approach 
presented in \cite{CaBKW} consists in resuming the idea of E.~Grenier,
writing the unknown phase $\Phi^\eps$
as
\begin{equation*}
  \Phi^\eps = \phi_{\rm eik} + \phi^\eps,
\end{equation*}
and considering \eqref{eq:systemmanuel} where $\nabla \phi^\eps$ has replaced
$\nabla \Phi^\eps$ as an unknown function. This procedure is similar
to linearizing \eqref{eq:systemmanuel} in $\Phi^\eps$, around
$\phi_{\rm eik}$. Of course, extra terms appear at this stage. Note
that the space where we seek $\Phi^\eps$ is of mixed type: $\Phi^\eps$
is the sum of a smooth, sub-quadratic (and possibly unbounded)
function, and the phase $\phi^\eps(t,\cdot)$ which is sought in Sobolev spaces
$H^s(\R^n)$. Nevertheless, $\phi^\eps$ must not be considered as
small, as shown by the analysis of \S\ref{sec:grenier}. 
\smallbreak

The good news is that the extra terms that have appeared can be treated
as \emph{semi-linear} perturbations in the energy estimates. This is
due to the fact that the 
phase $\phi_{\rm eik}$ is sub-quadratic in space. Therefore, the
analysis of \S\ref{sec:grenier} is easily adapted: provided that we
assume $f'>0$, an analogue of Theorem~\ref{theo:emmanuel} is
available. Note that unless $f'=\text{Const.}$ (in which case the
symmetrizer $S$ is constant), we need the extra decay assumption on
the initial profile:
  \begin{equation*}
xa_0\in \cap_{s\ge 0}H^s,\quad \text{and }    
    xa_0^\varepsilon \to xa_0 \text{ in }H^s(\R^n),\quad \forall s\ge
    0. 
  \end{equation*}
In particular, a local solution to \eqref{eq:eulerV} is constructed. 
We refer to \cite{CaBKW} for precise statements in this case. 
\begin{remark}
  For Schr\"odinger--Poisson equations in space dimension
  $n\ge 3$, the idea of E.~Grenier was adapted in \cite{AC-SP}, under
  more general geometrical assumptions. For instance, 
  solutions that do not necessarily have a zero limit at spatial
  infinity are considered. Under the assumptions of \cite{ZhangSIMA},
  a point-wise asymptotics of the wave function is given, which is
  more precise that the results in  \cite{ZhangSIMA}. 
\end{remark}

\subsection{The case of focusing nonlinearities}
\label{sec:foc}
Note that in \S\ref{sec:grenier}, the study of \eqref{eq:systemmanuel}
involves a quasi-linear system whose principal part writes:
\begin{equation*}
  \square_{f'}=
\d_{t}^{2}-\DIV\(f'(\left\lvert u^\eps\right\rvert
^2)\nabla\cdot\). 
\end{equation*}
This has the same form as the principal part for
\eqref{eq:systemmanuellimite}, which is the limiting system expected
in general, whichever formal approach is followed. When $f'>0$, we
face a quasi-linear wave equation. We have pointed out some open
problems under the weaker assumption $f'\ge 0$ (a case where loss of
hyperbolicity may occur). When $f'<0$, the above
operator becomes elliptic: it does not seem 
adapted to work in Sobolev spaces any more. On the other hand, data
and solutions with analytic regularity seem appropriate. 
\smallbreak

In \cite{PGX93}, P.~G\'erard works with the analytic regularity, when
the space variable $x$ belongs to the torus $\T^n$, without external
potential ($V\equiv 0$). Note
that the only assumption needed on the nonlinearity $f$ is analyticity
near the range of $|a_0|^2$. This includes the focusing case $f'<0$,
as well as the defocusing quintic case $f(y)=y^2$ for instance.
\smallbreak

The initial phase $\phi_0$ is supposed real analytic, and the initial
amplitude is analytic in the sense of J.~Sj\"ostrand \cite{SjoAst}: there
exist $\ell>0$, $A>0$, $B>0$ such that, for all $j\ge 0$, $a_j$ is
holomorphic in $\{|\IM x|<\ell\}$, and 
\begin{equation*}
  |a_j(x)|\le A B^j j!
\end{equation*}
Denoting $\overline a(t,x)$ the complex conjugate of $a(\overline
t,\overline x)$, P.~G\'erard considers the system:
\begin{equation*}
  \left\{
    \begin{aligned}
      \d_t v^\eps &=-v^\eps\cdot \nabla v^\eps -\nabla
      f\(a_0\overline a_0\),\\ 
\d_t a^\eps &=  -v^\eps\cdot \nabla a^\eps-\frac{1}{2}a^\eps \DIV
v^\eps + i\frac{\eps}{2}\Delta a^\eps -\frac{ia^\eps }{\eps}\(
f\(a^\eps \overline a^\eps\)- f\(a_0 \overline a_0 \)\). 
    \end{aligned}
\right.
\end{equation*}
A solution of the form 
\begin{equation*}
 u^\eps = a^\eps e^{i\phi/\eps},\quad a^\eps(t,x) =\sum_{j\ge 0}
\eps^j a^{(j)} (t,x) ,
\end{equation*}
where the sum is defined in the sense of J.~Sj\"ostrand, is thus
obtained. Setting 
\begin{equation*}
  v^\eps =  e^{i\phi/\eps}\sum_{j\le 1/(C_0 \eps)}\eps^j a^{(j)}
\end{equation*}
for $C_0$ sufficiently large, the approximate solution $v^\eps$
satisfies:
\begin{equation*}
  i\eps \d_t v^\eps +\frac{\eps^2}{2}\Delta v^\eps =
  f\(|v^\eps|^2\)v^\eps+\O\(e^{-\delta /\eps}\), 
\end{equation*}
for some $\delta >0$. Essentially, this source term is sufficiently
small to overcome the difficulty pointed out at the beginning of
\S\ref{sec:strong}: for small time independent of $\eps$, the
exponential growth provided by 
Gronwall lemma is more than compensated by the term $e^{-\delta
  /\eps}$. We refer to \cite{PGX93} for precise statements and
elements of proof. 
\subsection{The integrable case}
\label{sec:int}
In the one-dimensional case, $n=1$, for a cubic nonlinearity
($f(y)=\pm y$), the Schr\"odinger equation is completely integrable. This
property remains with a time-independent external potential which is a
polynomial of degree at most two  \cite[p.~375]{AblowitzClarkson}. 
\smallbreak

In the absence of potential, several papers have studied the
semi-classical limit for \eqref{eq:nlssemi}--\eqref{eq:CIbkw} for the
cubic NLS in space dimension one. See for instance 
\cite{JLM} in the defocusing case, and
\cite{KMM,TVZ} in the focusing case. A very interesting aspect of this
approach is that it yields a description of the solution $u^\eps$ even
after shocks have appeared for the limiting Euler equation
\eqref{eq:euler}. This description involves theta functions, and the
so-called Whitham equations (see \cite{TianYe}). In particular, this
approach seems to confirm the formal discussion of the end of
\S\ref{sec:grenier}: in the defocusing case, the $L^\infty$ norm of
the solution $u^\eps$ 
remains bounded as $\eps\to 0$, for all time.  
\smallbreak

Unfortunately, it seems that all the results in the integrable case
have been written in a way that makes any comparison with the other
results mentioned above very difficult. The last step of inverse
scattering is not always performed, which should yield a point-wise
asymptotics of the wave function $u^\eps$. Moreover, the spaces in
which it would be available are not completely clear. The space
$L^\infty_{\rm loc}(\R_x)$ seems the most natural candidate. A bridge
between the approaches of \S\ref{sec:brenier} and \S\ref{sec:grenier} on
the one hand, and the approaches in the integrable case on the other
hand, would certainly be welcome in the community of
semi-classical analysis for nonlinear Schr\"odinger equations.

\section{Propagation of concentrated initial data}
\label{sec:concen}

\subsection{Defocusing nonlinearity}
\label{sec:concdefoc}
We now consider \eqref{eq:nlssemiconc} with $\l>0$. By scaling, we may
assume $\l=1$. The general heuristic argument is the following. For
$t$ close to zero, the solution $\u^\eps$ remains concentrated near
the point $x_0$, at a scale of order $\eps$. Since the potential $V$
does not depend on $\eps$, we have $V \u^\eps \sim V(x_0)\u^\eps$: the
potential can be considered as constant at leading order. Introduce
the function $\psi^\eps$ given by the scaling
\begin{equation*}
  \u^\eps(t,x) = \frac{1}{\eps^{n/2}}\psi^\eps\( \frac{t}{\eps}\virgp
  \frac{x-x_0}{\eps}\)e^{i\(x_0\cdot \xi_0/\eps -
  V(x_0)t/\eps \)}.  
\end{equation*}
The Cauchy problem \eqref{eq:nlssemiconc} is equivalent to:
\begin{equation*}
  i\d_t \psi^\eps +\frac{1}{2}\Delta \psi^\eps = \(V(x_0+\eps
  x)-V(x_0)\)\psi^\eps +|\psi^\eps|^{2\si}\psi^\eps \quad ;\quad
  \psi^\eps(0,x) = R(x)e^{ix\cdot \xi_0}. 
\end{equation*}
The above argument suggests that we have $\psi^\eps \sim \psi$, where
$\psi$ is independent of $\eps$ and solves:
\begin{equation}
  \label{eq:psi}
  i\d_t \psi +\frac{1}{2}\Delta \psi = |\psi|^{2\si}\psi \quad ;\quad
  \psi(0,x) = R(x)e^{ix\cdot \xi_0}. 
\end{equation}
Under suitable assumptions on $\si$ and $R$, there is scattering for this
equation (see e.g.\  \cite{CazCourant,CW92,NakanishiOzawa}): there
exist $\psi_\pm \in L^2(\R^n)$ such that
\begin{equation}\label{eq:scatt}
  \left\| \psi(t) - e^{i\frac{t}{2}\Delta} \psi_\pm\right\|_{L^2}\to 0 \quad
  \text{as }t\to \pm \infty. 
\end{equation}
The standard asymptotics of the free Schr\"odinger group
$e^{i\frac{t}{2}\Delta}$ then yields:
\begin{equation*}
  \psi(t,x)\Eq t {\pm \infty}
  \frac{e^{i|x|^2/(2t)}}{(it)^{n/2}}\widehat \psi_\pm \(\frac{x}{t}\),
\end{equation*}
where the Fourier transform is given by
\begin{equation*}
  \F f(\xi)=\widehat f(\xi)=\frac{1}{(2\pi)^{n/2}}\int_{\R^n}e^{-ix\cdot
  \xi}f(x)dx. 
\end{equation*}
Back to $\u^\eps$, this yields, for $t\gg \eps$ and so long as we consider
the external 
potential as constant: 
\begin{equation}\label{eq:14h50}
  \u^\eps(t,x) \sim \frac{1}{(it)^{n/2}}\widehat\psi_+\(
  \frac{x-x_0}{t}\)e^{i\frac{|x-x_0|^2}{2\eps t}}
e^{i\(x_0\cdot \xi_0/\eps - 
  V(x_0)t/\eps \)}.
\end{equation}
Indeed, we have the following rigorous result: 
\begin{proposition}[\cite{CM}, Proposition~6.3]\label{prop:1}
Let $V$ satisfying Assumption~\ref{hyp:geom}. Let
$R\in \Sigma:=H^1 \cap \F(H^1)$, and 
\begin{equation*}
  \frac{2-n +\sqrt{n^2+12n+4}}{4n}\le \si<\frac{2}{n-2}\cdot
\end{equation*}
Then for any $\Lambda >0$, the following holds:\\
$1.$ There exists $\eps(\Lambda)>0$ such that for $0<\eps\leq
   \eps(\Lambda)$, the initial value problem
\eqref{eq:nlssemiconc} 
has a unique solution $\u^\eps \in
   C([-\Lambda\eps,\Lambda\eps];\Sigma)$. \\
$2.$ This solution satisfies the following asymptotics,
\begin{equation*}
\limsup_{\eps \to 0}\sup_{|t|\le \Lambda \eps}
\left\|\u^\eps(t)-\v^\eps(t) \right\|_{L^2}=0\ ,
\end{equation*}
\begin{equation*}
\text{where
}\v^\eps\text{ is given by}\quad \v^\eps(t,x)=\frac{1}{\eps^{n/2}}\psi
\left(\frac{t}{\eps}\virgp\frac{x-x_0}{\eps} \right)e^{i\(x_0\cdot
  \xi_0/\eps -  
  V(x_0)t/\eps \)},   
\end{equation*}
and $\psi\in C(\R;\Sigma)$ is given by \eqref{eq:psi}. 
\end{proposition}
A transition is expected to occur in the above boundary layer, that is
for $|t|=\Lambda \eps$ and $\Lambda \gg 1$. The heuristic argument
consists in saying that because of dispersion for $\psi$, the external
potential $V$ can no longer be considered as constant. On the other
hand, and for the same reason, the nonlinearity ceases to be relevant
at leading order: for $\eps\ll \pm t\le T$, we expect $\u^\eps \sim
\u_\pm^\eps$, where 
\begin{equation}
i\eps \d_t \u_\pm^\eps +\frac{\eps^2}{2}\Delta \u_\pm^\eps =
V\u_\pm^\eps\quad ;\quad 
\u^\eps(0,x) = \frac{1}{\eps^{n/2}}\psi_\pm
\left( 
\frac{x}{\eps}\right)e^{i\(x_0\cdot
  \xi_0/\eps -  
  V(x_0)t/\eps \)},
\end{equation}
and $\psi_\pm$ are given by \eqref{eq:scatt}. The value of $T$ is not
arbitrary: the asymptotic behavior of $\u_\pm^\eps$ involves the
classical trajectories associated to $V$. These trajectories may
refocus at one point; this is the case when $V$ is an isotropic
harmonic potential for instance. 
\smallbreak

Proving the above asymptotics for $\eps\ll \pm t\le T$ is actually an
open problem for general potentials satisfying
Assumption~\ref{hyp:geom}, even for time-independent potentials. It
has been proved when $V=V(x)$ is exactly a polynomial of degree at
most two, in \cite{CaIHP} for the case of refocusing(s), and in  \cite{CM}
for the complementary case. 
\smallbreak
The restriction to this class of polynomial potentials is certainly
purely technical, and we know explain it. The proof of the asymptotics
for $\eps\ll \pm t\le T$ relies on the use of operators well suited to
the propagation of classical trajectories associated to $V$. In the
linear setting, good candidates to meet this requirement are given by
the action of Heisenberg derivatives (see e.g.\  \cite{Robert}):
\begin{equation*}
  U^\eps(t)\eps\nabla U^\eps(-t)\text{
  and } U^\eps(t)\frac{x-x_0}{\eps}U^\eps(-t), \quad \text{where }U^\eps(t) =
  e^{-i\frac{t}{\eps}\(-\frac{\eps^2}{2}\Delta +V\)}. 
\end{equation*}
The main technical remark in \cite{CaIHP,CM} is that when $V$ is a
polynomial of degree at most two, then the above two Heisenberg
derivatives are very interesting for nonlinear problems too. Indeed,
we can find $p=p(t)$, and $\phi=\phi(t,x)$ real-valued, such that, for
instance:
\begin{equation}\label{eq:heisen}
  U^\eps(t)\frac{x-x_0}{\eps}U^\eps(-t) =
  p(t)e^{i\phi(t,x)/\eps}\nabla\( e^{-i\phi(t,x)/\eps} \cdot\). 
\end{equation}
In \cite{CM}, it is proved that an operator of the form of the 
right hand side of \eqref{eq:heisen} commutes with $U^\eps(t)$ if and
only if $V$ is a 
polynomial of degree at most two, and $\phi$ solves the eikonal
equation \eqref{eq:eik}. The fact that an Heisenberg derivative
commutes with the group $U^\eps(t)$ is a straightforward consequence
of its definition. The right hand side of \eqref{eq:heisen} implies
two important things:
\begin{itemize}
\item This Heisenberg derivative acts on gauge invariant
  nonlinearities $G(|u|^2)u$ like a derivative.
\item Weighted Gagliardo--Nirenberg inequalities are available, of the
  form 
  \begin{equation*}
    \|\varphi\|_{L^r}\le C_r
    |p(t)|^{-\delta(r)}\|\varphi\|_{L^2}^{1-\delta(r)}\left\|
    U^\eps(t)\frac{x-x_0}{\eps}U^\eps(-t)\varphi \right\|_{L^2}^{\delta(r)}. 
  \end{equation*}
\end{itemize}
To illustrate the use of these properties, we recall
\cite[Corollary~1.3]{CaIHP}: 
\begin{proposition}\label{prop:harmo}
Let $R\in \Sigma$. 
Assume that $\u^\eps$ solves \eqref{eq:nlssemiconc} with
$x_0=\xi_0=0$ and 
\begin{equation*}
  V(x)=\frac{|x|^2}{2}\cdot
\end{equation*}  
Let $\psi_\pm =W_\pm^{-1}R$ be given by \eqref{eq:scatt} (upon
suitable assumptions on $\si$, see e.g.\  Prop.~\ref{prop:1}). 
Then for any $2<r<\frac{2n}{n-2}$, the following asymptotics 
holds in $L^2\cap L^r$:  
\begin{itemize}
\item If $0< t<\pi$, then 
$\displaystyle \u^\eps(t,x) \Eq \eps 0 \left(\frac{1}{i \sin
t}\right)^{n/2}\widehat{\psi_+}\left(\frac{x}{\sin t} \right)
e^{i\frac{|x|^2}{2\eps\tan t}}.$
\item If $-\pi< t<0$, then 
$\displaystyle \u^\eps(t,x) \Eq \eps 0 \left(\frac{1}{i \sin
t}\right)^{n/2}\widehat{\psi_-}\left(\frac{x}{\sin t} \right)
e^{i\frac{|x|^2}{2\eps\tan t}}.$
\end{itemize}
\end{proposition}
\begin{remark}
The result of \cite{CaWigner} shows that in the above case, Wigner
measure is not a good tool to characterize the behavior of
$\u^\eps$. More precisely, we can find $R_1,R_2\in \Sigma$ such that
the Wigner measures for the corresponding solutions $\u_1^\eps$ and
$\u_2^\eps$ coincide at time $t=-\pi/2$, but are different at time
$t=\pi/2$. The crossing of a focal point may lead to an ill-posed
Cauchy problem as far as Wigner measures are concerned. 
\end{remark}
We see that the formal asymptotics \eqref{eq:14h50} is valid only in
the transition r\'egime $t=\Lambda \eps$, with $\Lambda \gg 1$. For
larger times, the trigonometric functions in the above result account
for the dynamical influence of the harmonic potential. 
\smallbreak

In the above case of an isotropic harmonic potential, the above result
can be iterated in time. Recall that the (nonlinear) scattering
operator $S$ associated to \eqref{eq:psi} maps $\psi_-$ to $\psi_+$,
given by \eqref{eq:scatt}. 
\begin{corollary}
  Under the assumptions of Proposition~\ref{prop:harmo}, consider $k\in
  \N$. For $k\pi<t<(k+1)\pi$, and $2<r<\frac{2n}{n-2}$, the following
  asymptotics  
holds in $L^2\cap L^r$: 
\begin{equation*}
  \u^\eps(t,x) \Eq \eps 0 \frac{e^{-in\frac{\pi}{4}-ink\frac{\pi}{2}}}{|\sin
t|^{n/2}}\widehat{S^k\psi_+}\left(\frac{x}{\sin t} \right)
e^{i\frac{|x|^2}{2\eps\tan t}},
\end{equation*}
where $S^k$ denotes the $k^{\rm th}$ iterate of the scattering
operator $S$. 
\end{corollary}
The phase shift $e^{-ink\frac{\pi}{2}}$ corresponds to successive
Maslov indices: this is a linear phenomenon \cite{Du}. On the other
hand, we see 
that a nonlinear phenomenon occurs at leading order at time $t=k\pi$,
which is measured by the scattering operator $S$. 

\subsection{Focusing nonlinearity}
\label{sec:concfoc}

When $\l<0$ in \eqref{eq:nlssemi}--\eqref{eq:CIconc}, we assume
similarly that $\l=-1$. We let $R=Q$, the unique
positive, radially symmetric (\cite{Kwong}), solution of: 
\begin{equation*}
  -\frac{1}{2}\Delta Q +Q = |Q|^{2\si}Q. 
\end{equation*}
Now, the focusing nonlinearity is an obstruction to dispersive
phenomena. The
solution $u^\eps$ is expected to keep the ground state as a leading
order profile. Nevertheless, the point where it is centered in the
phase space, initially $(x_0,\xi_0)$, should evolve according to the
Hamiltonian flow \eqref{eq:hamilton}.  In the absence of external
potential, $V\equiv 0$, we have explicitly:
\begin{equation*}
  u^\eps(t,x) = Q\( \frac{x-x(t)}{\eps}\) e^{ix\cdot \xi(t)/\eps +
  i\theta(t)/\eps}, 
\end{equation*}
where $(x(t),\xi(t))=(x_0 +t\xi_0,\xi_0)$ solves \eqref{eq:hamilton}
with initial data 
$(x_0,\xi_0)$, and $\theta(t) = t- t|\xi_0|^2/2$. When $V$ is not
trivial, seek $u^\eps$ of the form of a rescaled WKB expansion:
\begin{equation*}
  u^\eps(t,x)\sim \(\sum_{j\ge 0} \eps^jU_j\(\frac{t}{\eps}\virgp
  \frac{x-x(t)}{\eps}  
  \)\)e^{i\phi(t,x)/\eps}. 
\end{equation*}
Note that this scaling meets the exact result of the case $V\equiv
0$. Plugging this expansion into \eqref{eq:nlssemi}--\eqref{eq:CIconc} and
canceling the $\O(\eps^0)$ term, we get:
\begin{equation*}
  i\d_t U_0 +\frac{1}{2}\Delta U_0 + U_0\( -\d_t \phi
  -\frac{1}{2}|\nabla \phi|^2 -V +|U_0|^{2\si}\) -i \(\dot x(t)-\nabla
  \phi\)\cdot \nabla U_0=0. 
\end{equation*}
Impose the leading order profile to be the standing wave given by
\begin{equation*}
  U_0(t,x) = e^{it}Q(x). 
\end{equation*}
Then the above equation becomes:
\begin{equation*}
  U_0\( -\d_t \phi
  -\frac{1}{2}|\nabla \phi|^2 -V \) -i \(\dot x(t)-\nabla
  \phi\)\cdot \nabla U_0=0. 
\end{equation*}
Since $U_0 e^{-it}$ is real-valued, and since we seek a real-valued phase
$\phi$, this yields:
\begin{align*}
  &\d_t \phi
  +\frac{1}{2}|\nabla \phi|^2 +V =0\quad ;\quad \phi(0,x)=x\cdot
  \xi_0.\\
& \dot x(t)=\nabla
  \phi(t,x). 
\end{align*}
The first equation is the eikonal equation
\eqref{eq:eik}. We infer that we have exactly 
\begin{equation*}
  \nabla \phi\( t,x(t)\) = \xi(t). 
\end{equation*}
The form of $U_0$ and the exponential decay of $Q$ show that we can
formally assume that $x=x(t) +\O(\eps)$. In this case, 
\begin{equation*}
  \nabla
  \phi(t,x)=\nabla
  \phi\(t,x(t)\)+ \O(\eps) = \xi(t) + \O(\eps) = \dot x(t)+ \O(\eps) .
\end{equation*}
Thus, we have canceled the $\O(\eps^0)$ term, up to adding extra terms
of order $\eps$, that would be considered in the next step of the
analysis, which we stop here. Back to $u^\eps$, this formal
computation yields
\begin{equation*}
  u^\eps(t,x) \sim Q\(\frac{x-x(t)}{\eps}\)e^{i\phi(t,x)}\sim
  Q\(\frac{x-x(t)}{\eps}\) e^{ix\cdot \xi(t)/\eps +i\theta(t)/\eps},
\end{equation*}
where $\theta(t) = \dis t\(1- |\xi_0|^2/2-V(x_0)\) +\int_0^t x(s)\cdot
\nabla V(x(s))ds$. 
\smallbreak

To give the above formal analysis a rigorous justification, the
following assumptions are made in \cite{KeraaniAA}:
\begin{hyp}\label{hyp:sahbi}
  The nonlinearity is $L^2$-subcritical: $\si<2/n$.\\
The potential $V=V(x)$ is real-valued, and can be written as
$V=V_1+V_2$, where
\begin{itemize}
\item $V_1\in W^{3,\infty}(\R^n)$.
\item $\d^\alpha V_2\in W^{2,\infty}(\R^n)$ for every multi-index
  $\alpha$ with $|\alpha |=2$.  
\end{itemize}
\end{hyp}
For instance, $V$ can be an harmonic potential. 
\begin{theorem}[\cite{KeraaniAA}]\label{theo:sahbi}
  Let $x_0,\xi_0\in \R^n$. Under Assumption~\ref{hyp:sahbi}, the
  solution $u^\eps$ to 
  \eqref{eq:nlssemi}--\eqref{eq:CIconc} with $R=Q$ can be approximated
  as follows:
  \begin{equation*}
    u^\eps(t,x) = Q\(\frac{x-x(t)}{\eps}\) e^{ix\cdot \xi(t)/\eps + i
    \theta^\eps(t)/\eps} +\O(\eps) \quad \text{in }L^\infty_{\rm
    loc}(\R_t;X^\eps), 
  \end{equation*}
where $(x(t),\xi(t))$ is given by the Hamiltonian flow, the
real-valued 
function $\theta^\eps$ depends on $t$ only, and $X^\eps$ is defined by
the norm
\begin{equation*}
  \|f\|_{X^\eps}^2=\frac{1}{\eps^n}\|f\|_{L^2}^2
  +\frac{1}{\eps^{n-2}}\|\nabla f\|_{L^2}^2  . 
\end{equation*}
\end{theorem}
\begin{remark}
  The assumption $\si<2/n$ is crucial for the above result to
  hold. Indeed, if $\si=2/n$ and $V$ is the isotropic harmonic
  potential 
  \begin{equation*}
    V(x)=\frac{|x|^2}{2}\virgp
  \end{equation*}
then we have explicitly, when $x_0=\xi_0=0$ (see
\cite{CaM3AS,KeraaniAA}): 
\begin{equation*}
  u^\eps(t,x)= \frac{1}{(\cos t)^{n/2}} Q\( \frac{x}{\eps \cos t}\)
  e^{i\frac{\tan t}{\eps}-i\frac{|x|^2}{2\eps}\tan t},\quad 0\le
  t<\frac{\pi}{2}\virgp 
\end{equation*}
so the profile $Q$ is modulated as time evolves, in a fashion
similar to \S\ref{sec:concdefoc}. 
\end{remark}
The proof of the above result heavily relies on the orbital stability
of the ground state, which holds when $\si<2/n$. For $v\in H^1(\R^n)$,
denote
\begin{equation*}
  {\mathcal E}(v)= \frac{1}{2}\|\nabla v\|_{L^2}^2 -\frac{1}{\si
  +1}\|v\|_{L^{2\si+2}}^{2\si+2}.  
\end{equation*}
The ground state $Q$ is the unique solution, up to translation and
rotation,  to the minimization
problem:
\begin{equation*}
  {\mathcal E}(Q)= \inf \{ {\mathcal E}(v)\ ;\ v\in H^1(\R^n)\text{
  and }\|v\|_{L^2}=\|Q\|_{L^2}\}. 
\end{equation*}
The orbital stability is given by the following result:
\begin{proposition}[\cite{Weinstein85}]\label{prop:stab}
  Let $\si<2/n$. There exist $C,h>0$ such that if $\phi\in H^1(\R^n)$
  is such that $\|\phi\|_{L^2}=\|Q\|_{L^2}$ and ${\mathcal E}(\phi)-
  {\mathcal E}(Q) <h$, then:
  \begin{equation*}
    \inf_{y\in\R^n,\theta\in \T}\left\| \phi - e^{i\theta} Q(\cdot
    -y)\right\|_{H^1}^2 \le C \({\mathcal E}(\phi)-  {\mathcal E}(Q)\).
  \end{equation*}
\end{proposition}
The strategy in \cite{KeraaniAA} consists in applying the above result
to the function
\begin{equation*}
  v^\eps(t,x) = u^\eps\(t,\eps x+x(t) \)e^{-i(\eps x +x(t))\cdot
  \xi(t)/\eps}.  
\end{equation*}
For $A>0$ sufficiently large, let $\chi$ be a smooth non-negative
cut-off function, supported in $\{x\in \R^n; |x|\le 2A\}$, and
constant equal to $1$ in $\{x\in \R^n; |x|\le A\}$. Introduce the
error estimate 
$\eta^\eps(t)$ given by $\eta^\eps =
\eta^\eps_1+\eta^\eps_2+\eta^\eps_3+\eta^\eps_4$, where:
\begin{align*}
 &\eta^\eps_1(t)=\int_{\R^n} x\chi(x)m^\eps(t,x)dx- \|Q\|_{L^2}^2
 x(t),\\ 
 &\eta^\eps_2(t)= \int_{\R^n} \nabla V_2(x) m^\eps(t,x)dx -
 \|Q\|_{L^2}^2 \nabla V_2\( x(t)\),\\
&\eta^\eps_3(t)=\int_{\R^n} \xi^\eps(t,x)dx- \|Q\|_{L^2}^2
 \xi(t),\\
&\eta^\eps_4(t)=\int_{\R^n} \chi(x)V(x)m^\eps(t,x)dx- \|Q\|_{L^2}^2
 V\(x(t)\),\\
&m^\eps(t,x)= \frac{1}{\eps^n}|u^\eps(t,x)|^2\quad ;\quad \xi^\eps(t,x)
 =\frac{1}{\eps^{n-1}} \IM\( \overline u^\eps \nabla u^\eps\). 
\end{align*}
Noting that $\eta^\eps(0)=\O(\eps^2)$, the proof in \cite{KeraaniAA}
shows that $\eta^\eps(t)=\O(\eps^2)$ for $t\in [0,T_0]$ for some
$T_0>0$ independent of $\eps$. The proof eventually relies on
Gronwall lemma and a continuity argument. In order to invoke these
arguments, S.~Keraani uses Proposition~\ref{prop:stab} and the scheme
of the proof of J.~Bronski and R.~Jerrard \cite{BJ00}, based on
duality arguments and estimates on measures. Finally, the time $T_0$
given by the proof depends only on constants of the motion, so the
argument can be repeated indefinitely, to get the $L^\infty_{\rm loc}$
estimate of Theorem~\ref{theo:sahbi}. 
\smallbreak

In the particular case where the external potential $V$ is an
harmonic potential (isotropic or anisotropic), the proof can be simplified. We
invite the reader to pay attention to the short note
\cite{KeraaniCRAS}, where this simplification is available.  
\smallbreak

The phase shift $\theta^\eps$ in Theorem~\ref{theo:sahbi} is
not known in general. It is easy to guess from the arguments given
above that in the proof given by S.~Keraani, it stems from the use of
Proposition~\ref{prop:stab}. On the other hand, as noted in
\cite{KeraaniAA}, a time-dependent phase shift does not alter the
Wigner measure of $u^\eps$, which is an important physical quantity. 

\noindent {\bf Acknowledgments.} The author is grateful to
Thomas Alazard for his careful reading of the manuscript, and for his
comments.

\providecommand{\bysame}{\leavevmode\hbox to3em{\hrulefill}\thinspace}
\providecommand{\href}[2]{#2}

\end{document}